\documentclass[11pt]{article}
\usepackage[utf8]{inputenc}
\usepackage[margin=1in,left=1in]{geometry}
\usepackage{setspace}
\usepackage{graphicx}
\usepackage{amssymb, amsmath, amsthm, graphics}
\usepackage{mathabx,epsfig}
\usepackage[mathscr]{euscript}
\usepackage{tikz}
\usetikzlibrary{calc}
\usetikzlibrary{matrix}
\usepackage{array}
\newcolumntype{C}[1]{>{\centering\arraybackslash}p{#1}}

\usepackage{latexsym}
\usepackage[all]{xy}
\usepackage{color}

\usepackage{bbold}

\usepackage{tikz}
\input xy
\xyoption{all}
\usepackage{euscript}
\usepackage{multirow}
\usepackage{etex, pictexwd,dcpic}
\usetikzlibrary{positioning}
\usetikzlibrary{shapes.geometric}
\usetikzlibrary{shapes.misc}
\usetikzlibrary{calc}
\usetikzlibrary{positioning}

\usepackage{pgflibraryarrows}
\usepackage{pgflibrarysnakes}

\usetikzlibrary{trees} 
\usetikzlibrary[trees] 

\def\acts{\mathrel{\reflectbox{$\righttoleftarrow$}}}

\newtheorem{theorem}{Theorem}
\newtheorem{definition}[theorem]{Definition}

\newtheorem{corollary}[theorem]{Corollary}
\newtheorem{proposition}[theorem]{Proposition}

\newtheorem{remark}{Remark}

\newtheorem{example}{Example}

\numberwithin{equation}{section}

\renewcommand{\(}{\begin{equation*}}
\renewcommand{\)}{\end{equation*}}
\newcommand{\bea}{\begin{eqnarray*}}
\newcommand{\eea}{\end{eqnarray*}}

\def\endofproof {\hfill{$\Box$}\\}

\newcommand{\beq}{\begin{equation}}
\newcommand{\eeq}{\end{equation}}

\newcommand{\onto}{\twoheadrightarrow}

\newcommand{\into}{\hookrightarrow}

\newcommand{\theproof}{\noindent {\bf Proof.\ }}

\numberwithin{equation}{section}

\renewcommand{\(}{\begin{equation}}
\renewcommand{\)}{\end{equation}}

\def\1{{\bf 1}}

\def\<{\langle}
\def\>{\rangle}

\numberwithin{equation}{section}


\newcommand{\RR}{\ensuremath{\mathbb R}}

\newcommand{\HH}{\ensuremath{\mathbb H}}
\newcommand{\BB}{\ensuremath{\mathbf B}}

\newcommand{\sh}{\ensuremath{\mathscr{S}\mathrm{h}}}

\newcommand{\cartsp}{\mathscr{C}\mathrm{art}\mathscr{S}\mathrm{p}}

\newcommand{\map}{\mathrm{Map}}

\title{Cobordisms of global quotient orbifolds and an equivariant Pontrjagin-Thom construction}
\author{Daniel Grady}

\begin{document}

\maketitle
\begin{abstract}
We introduce an equivariant Pontrjagin-Thom construction which identifies equivariant cohomotopy classes with certain fixed point bordism classes. This provides a concrete geometric model for equivariant cohomotopy which works for any compact Lie group G. In the special case when G is finite or a torus, we show that our construction recovers the construction of Wasserman, providing a new perspective on equivariant bordism. We connect the results with bordisms of global quotient orbifolds, utilizing the machinery of Gepner-Henriques to describe bordisms of framed orbifolds in terms of equivariant cohomotopy.  We also illustrate the utility of the theory by applying our results to M-theory, thus connecting with recent work of Huerta, Sati and Schreiber.
\end{abstract}

\tableofcontents

\section{Introduction}

The idea of studying differentiable manifolds by means of homotopy theory goes back to the work of Lev Pontrjagin on mappings between spheres \cite{Po}. Almost 20 years after this initial work, both R\'ene Thom and Lev Pontrjagin realized that bordisms of manifolds could be completely classified using Pontrjagin's construction (Thom in the oriented case \cite{Th}, and Pontrjagin in the framed case \cite{Po2}). The result of Pontrjagin's work was a beautiful correspondence between the degree $k$ cohomotopy groups of a smooth manifold $M$ and the bordism classes of framed codimension-$k$ submanifolds of $M$. Since its conception, this correspondence has been generalized in a variety of ways, each of which retain the conceptual elegance of Pontrjagin's initial work. In the stable case \cite{Th}, the correspondence gives a geometric meaning to the homotopy groups of spheres and allows one to study the geometric theory of bordisms by means of homotopy theory. The categorification of the correspondence \cite{GMTW} can be used to provide fascinating equivalences between cohomology theories \footnote{For example, the $d=1$ dimensional case can be used to prove the main statement in \cite{BPQ}.} and has a number of applications in physics, in particular in TQFT's \cite{FH} \cite{At}. The geometric refinement of \cite{GS} provides a correspondence which allows for the addition of geometric data and has a wealth of applications to physics, in particular (non topological) QFT's. Finally, the correspondence gives rise to a complete classification of fully extended TQFT's  when generalized beyond spaces to $(\infty,n)$-categories.

Generalizing Thom's results to the equivariant case has shown a number of unforseen difficulties. In particular, the transversality theorems required to come up with a correspondence resembling the classical case do not always hold in the equivariant case (see \cite{Was}). Given the existence of an equivariant Whitney embedding theorem for $G$-manifolds \cite{Mos}\cite{Pal}, a crucial step in the Pontrjagin-Thom correspondence in the stable case, it is perhaps even more surprising that the classical correspondence does not carry through to the equivariant setting. These technicalities can be circumvented in a number of ways. For example, in \cite{Wan}, a Pontrjagin-Thom type theorem is established which related equivariant bordism to certain equivariant cohomology groups. However, in the framed case, the representing spectrum is not the equivariant sphere spectrum. On the other hand, it is also shown in \cite{Wan} that the equivariant sphere spectrum represents a modified version of the bordism category, which formally inverts the equivariant suspension maps, although this only holds in the case of finite groups. In a similar spirit, it is shown in \cite{BH} that after inverting certain tautological classes, equivariant bordism is represented by the corresponding equivariant spectrum. This answer, however, still does not address the question of what geometric theory the equivariant bordism spectrum represents in the general case of compact Lie groups, and this question remains open (see \cite[Section 6.2]{Sh} and \cite{Sh2} for discussion on this point).

The goal of this paper is twofold. The first goal is to address this question in the simple case of unstable, equivariant cohomotopy, akin to Pontrjagin's initial work. Surprisingly, it seems that there is very little literature which addresses this case (most are devoted to generalizations of the stable case), even though this is arguably the most crucial situation at a conceptual level. Historically, generalizations of the correspondence emerged from this case. The second goal is to recast the results in terms of orbifold cohomology for global quotient orbifolds. In \cite{GH} it is shown that the cohomology of orbifolds is intimately related to global equivariant homotopy theory via Elmendorf's theorem \cite{El}. Uniting these two perspective provides a power conceptual understanding of what happens in the equivariant case.

Finally, this work was inspired by the applications of cohomotopy in physics, a program initiated by Sati in \cite{Sa}, further developed in the setting of rational homotopy theory by Fiorenza, Sati and Schreiber in \cite{FSS1}, and developed in the rational equivariant setting with Huerta in \cite{HSS}. In the context of M-theory, M-branes arise as the fixed points of certain finite subgroups of $SU(2)$.  In \cite{FSS1}, it was shown that (rationally) these branes are classified by maps to the rational 4-sphere. As pointed pointed out in \cite{Sa}, there is good reason to believe that this should hold even non-rationally. Since cohomotopy classifies framed submanifolds, it is a tantilizing conjecture that perhaps this map to the 4-sphere singles out the framed bordism class of the branes. Even more, if one works out what happens equivariantly, one hopes that the branes will emerge as the fixed points of some action of a finite subgroup of $SU(2)$. In fact, our correspondence shows that indeed this is the case. It is striking that the correspondence takes this form, since in the initial stages (following the classical literature) it was not clear why the maps ought to classify fixed point submanifolds. 

The outline of the paper is as follows. In the first section we connect equivariant bordism with bordisms of global quotient orbifolds. The second section is the main section and is devoted to the construction of a collapse map which we call the fixed point Pontrjagin-Thom collapse map. This requires a structure on the normal bundle which is weaker than a framing, which we call a partial $V$-framing. We use this map to prove the following theorem.
\begin{theorem}
Let $M$ be a smooth compact $G$-manifold with $G$ a compact Lie group. Let $V$ be a $G$-representation with corresponding representation sphere $S^V$ and let $M/\!/G$ denote the resulting global quotient orbifold. Then for any closed normal subgroup $H\trianglelefteq G$, we have the following identifications
\begin{enumerate}
\item The fixed point Pontrjagin-Thom construction gives rise to an isomorphism
$$
{\rm PT}_{\rm fix}:\pi_H^V(M) \cong \left\{\begin{array}{c}
\text{Bordism classes of $d_H$-dimensional }
\\
\text{ $H$-fixed point submanifolds of $M$ }  
\\
\text{with partial $V$-framed normal bundles.}
\end{array}\right\}\;,
$$
where $d_H={\rm dim}(M^H)-{\rm dim}(V^H)$.
\item We have an isomorphism
$$\pi_0\vert {\bf Map}_{\BB G}(M/\!/H,S^V/\!/H)\vert \cong  \pi_H^V(M)\;,$$
where $\pi^V_H(M)$ is the unstable $H$-equivariant cohomotopy in ${\rm RO}(G)$-degree $V$, regarded as an $H$-representation, ${\bf Map}_{\BB G}(M/\!/H,S^V/\!/H)$ denotes the mapping stack in the slice $\infty$-topos $\sh_{\infty}(\cartsp)_{\BB G}$, and the bars $\vert \cdot \vert$ denote the geometric realization of this stack as a space.
\item There is a bijection 
$$
\bigoplus_{[H]\leq G}\left\{\begin{array}{c}
\text{Bordism classes of $d_H$-dimensional }
\\
\text{ $H$-fixed point submanifolds of $M$ }  
\\
\text{with partial $V$-framed normal bundles}
\end{array}\right\}\cong \left\{\begin{array}{c}
\text{Bordism classes of $d$-dimensional} 
\\
\text{fixed point suborbifolds of $M/\!/G$ with}
\\
\text{partial $V$-framed normal bundles.}
\end{array}\right\}
$$
where $H$ runs over all conjugacy classes of closed normal subgroups \footnote{Note that since $G$ is a compact Lie group, there are at most countably many conjugacy classes of closed subgroups.} of $G$, $d_H={\rm dim}(M^H)-{\rm dim}(V^H)$ and $d$ ranges over the dimensions $d_H$ as $H$-varies. 
\end{enumerate}
\end{theorem}
\noindent The last section discusses the particular case of $SU(2)$ acting on the representation sphere $S^{\mathbb{H}}$, with $\mathbb{H}$ the quaternions. We apply the result to various ADE actions on 11-dimensional space-time arising in $M$-theory, connecting with the work of \cite{HSS}, but in the general (i.e. non-rational) setting. 

\medskip Since we will be dealing with various $G$-actions, $G$-fixed points, $G$-invariant submanifolds, and equivariant maps which could easily be confused, we will standardize our notation and terminology as follows.
\begin{center}
\begin{tabular}{l*{2}{C{5cm}}C{10cm} }
\hline
{\bf Terminology} & {\bf Notation} & {\bf Meaning}
\\
\hline
$G$-manifold & $G\acts M$ or simply $M$ & $G$ acts on $M$ by diffeomorphisms
\\
\hline
$G$-invariant submanifold & $(G\acts X)\into (G\acts M)$ or $X\into M$& $X$ is a submanifold which is closed under the $G$-action
\\
\hline 
full fixed point submanifold of $M$ & $M^G$ & $M^G$ is the submanifold of $M$ whose underlying set is the fixed point set of $G$
\\
\hline
a fixed point submanifold of $M$ & $X^G$ or explicitly $X^G\into M^G\into M$ & $X^G$ is a submanifold of the full fixed point submanifold $M^G$
\\
\hline 
global quotient orbifold & $M/\!/G$ & $M/\!/G$ is the quotient stack by a $G$-action on $M$ 
\end{tabular}
\end{center}

\section{Orbifolds and cobordism}

Throughout, we will take the Lie groupoid perspective on orbifolds \cite{Mo},\cite{MP}. Thus, what we mean by an orbifold is a Lie groupoid which is equivalent to proper \'etale Lie groupoid. In particular, the propety \'etale forces the automorphism groups to be discrete. Properness implies in addition that this discrete group is finite. 

\medskip
In the ordinary case of smooth manifolds, the classical Pontrjagin-Thom construction gives a bijective correspondence between framed cobordism classes of codimension-$k$ submanifolds, embedded in some ambient manifold $M$, and the unstable cohomotopy set $\pi^k(M)$ \footnote{Note that this set does not admit a canonical group structure.}. Clearly, at the most basic level, we need to know what it means for one orbifold to embed into another. Surprisingly, a universally agreed upon definition of such a basic concept seems to be lacking in the literature. We offer a definition which is natural in the following sense. An embedding of smooth manifolds is a monomorphism in the category of smooth manifolds which is also an immersion and an embedding of underlying topological spaces \footnote{Monomorphisms in the category of smooth manifolds (with smooth maps as morphisms) are injective smooth maps. Hence, the extra condition of this injection being an immersion is necessary.}. Monomorphisms in the category of Lie groupoids are smooth functors which are injective on both objects an morphisms. Thus, in analogy with manifolds, we make the following natural definition.
\begin{definition}\label{emboorbis}
A morphism $i:\mathscr{X}\to \mathscr{Y}$ of orbifolds is called an embedding if the induced smooth maps
$$i_1:\mathscr{X}_1\to \mathscr{Y}_1,\ \ i_0:\mathscr{X}_0\to \mathscr{Y}_0\;,$$
defined on morphisms and objects, are embeddings of smooth manifolds.
\end{definition}

\begin{remark}
There is an intrinsic notion of embedding in any $\infty$-topos exhibiting differential cohesion \footnote{I am grateful to Urs Schreiber for pointing this out.}. In the $\infty$-category of sheaves on cartesian spaces (with infintesimal discs adjoined), one can show that any embedding of orbifolds whose underlying mono is presented by a 1-categorical mono satisfies precisely the conditions of definition \ref{emboorbis}. 
\end{remark}

A \emph{suborbifold} $\mathscr{X}\into \mathscr{Y}$ is an equivalence class of embeddings of $\mathscr{X}$ into $\mathscr{Y}$. Such an equivalence class admits a canonical representative -- namely, a pair of submanifolds $\mathscr{X}_1\subset \mathscr{Y}_1$ and $\mathscr{X}_0\subset \mathscr{Y}_0$ such that $s(\mathscr{X}_1),t(\mathscr{X}_1)\subset \mathscr{X}_0$, where $s$ and $t$ are the source and target maps for $\mathscr{Y}$, gives rise to parallel arrows $s,t:\mathscr{X}_1\to \mathscr{X}_0$ (by restriction). If the quadruple $\mathscr{X}=(\mathscr{X}_1,\mathscr{X}_0,s,t)$ admits the structure of a Lie groupoid, then it is a suborbifold. Moreover, it is easy to show that any embedding $i:\mathscr{X}\into \mathscr{Y}$ is canonically equivalent to an orbifold of this type (simply consider the induced image submanifolds of objects and morphisms). As a consequence of our definition \ref{emboorbis}, we have the following immediate examples of suborbifolds.

\begin{example}
Let $G$ be finite group acting smoothly on a compact manifold $M$. Then the quotient Lie groupoid $M/\!/G$ is an orbifold. The smooth manifold of morphisms is given by $M\times G$ and the objects are $M$. A suborbifold of this groupoid is a pair $X,H$ with $H$ a subgroup of $G$  and $X$ an $H$-invariant submanifold of $M$.
\end{example}

\begin{example}\label{gloqsub}
Let $G$ be a compact Lie group and let $M$ be a compact $G$-manifold. Suppose the action groupoid $M/\!/G$ has finite isotropy groups at each point $x\in M$. Then $M/\!/G$ is an orbifold. Every suborbifold of $M/\!/G$ is of the form $X/\!/H$ with $H$ a closed subgroup and $X$ an $H$-invariant submanifold.
\end{example}

\begin{remark}
We will use the notation $M/\!/G$ rather loosely. In some instances, we will take this to mean the groupoid associated with the action, while in others we will take this to mean the quotient stack (or $\infty$-stack). The former is a presentation for the later. If the need arises to distinguish between the flexible homotopy quotient and the strict 1-categorical model, we will say so. 
\end{remark}

The suborbifolds which arise in the global quotient case are thus themselves global quotient orbifolds. Moreover, the isotropy groups of $M/\!/G$ govern the local structure of these suborbifolds completely. Indeed, let $M/\!/G$ be a global quotient orbifold and let $X/\!/H$ be a suborbifold. Then each orbifold atlas $\{U_i,G_{x_i},\phi_i\}$ of $M/\!/G$ restricts to an atlas of $X/\!/H$ of the form $\{V_i,H_{x_i},\varphi_i\}$ with $V_i=U_i\cap N$, $H_{x_i}=G_{x_i}\cap H$ and $\varphi_i=\phi_i\vert_{X}$.

Bordisms of orbifolds have been studied before in the literature (although not extensively). In \cite{Ang} the notion of $(\mathcal{F},\mathcal{F}^{\prime})$-orbifold bordism is introduced, where $\mathcal{F}$ and $\mathcal{F}^{\prime}$ are families of local representations. In the case where $\mathcal{F}^{\prime}=\emptyset$, this theory reduces to the usual notion of bordism. In particular, two $n$-dimensional orbifolds $\mathscr{X}^n$ and $(\mathscr{X}^{\prime})^n$ are bordant if there is an $(n+1)$-dimensional orbifold $\mathscr{W}^{n+1}$ and an \emph{isomorphism} of orbifolds $\partial \mathscr{W}^{n+1}\cong \mathscr{X}^n\sqcup (\mathscr{X}^{\prime})^n$. In \cite{Joy}, the notion of bordism of orbifolds is weakened to fit naturally into the higher categorical context. In particular, one only requires that $\partial \mathscr{W}^{n+1}\simeq \mathscr{X}^n\sqcup (\mathscr{X}^{\prime})^n$ is an \emph{equivalence}. 

Since our orbifolds are already embedded into a fixed global quotient orbifold $M/\!/G$, we do not need to distinguish between equivalence and isomorphism. More precisely, an \emph{embedded} orbifold $\mathscr{X}\subset \mathscr{Y}$ is already identified with its corresponding suborbifold, either up to isomorphism or equivalence. Thus, in the present case, the following definition is quite natural.

\begin{definition}
Let $\mathscr{X}^n,(\mathscr{X}^{\prime})^n\subset \mathscr{M}$ be embedded suborbifolds. Then we say that $\mathscr{X}^n$ and $(\mathscr{X}^{\prime})^n$ are {\bf bordant} if there is a suborbifold $\mathscr{W}^{n+1}\subset \mathscr{M}\times [0,1]$ with $\partial \mathscr{W}^{n+1}\subset \mathscr{M}\times \{0,1\}\cong \mathscr{M}\sqcup \mathscr{M}$ and such that $\partial \mathscr{W}^{n+1}= \mathscr{X}^n\sqcup (\mathscr{X}^{\prime})^n$. 
\end{definition}

Since all suborbifolds of a global quotient $M/\!/G$ are of the form $X/\!/H$ for some closed subgroup $H\leq G$, we immediately have the following proposition.
\begin{proposition}
Consider again the global quotient $M/\!/G$, where $G$ is a compact Lie group. Two suborbifolds of the form $X^n/\!/H$ and $(X^{\prime})^n/\!/H$ are bordant if and only if there is a closed subgroup $K\leq G$ and a $K$-invariant submanifold $W^{n+1}\subset M$ such that $\partial (W^{n+1}/\!/K)=(X^n\sqcup (X^{\prime})^n)/\!/H $.
\end{proposition}

Note that the previous proposition implies that we can take our bordism to be of the form $W^{n+1}/\!/H$.
\begin{corollary}\label{equitorbi}
Two suborbifolds $X^n/\!/H$ and $(X^{\prime})^{n}/\!/H$ are bordant if and only if they are bordant via a suborbifold of the form $W^{n+1}/\!/H$.
\end{corollary}
\theproof
Suppose that $X^n/\!/H$ and $(X^{\prime})^{n}/\!/H$ are bordant with bordism of the form $W^{n+1}/\!/K$, with $H\leq K$ a closed subgroup. Restricting the $K$ action to $H$ gives a bordism of the form $W^{n+1}/\!/H$.
\endofproof

Following the general theory of $\infty$-bundles (see \cite{NSS} for details) we have the following general definition
\begin{definition}
In any $\infty$-topos, a $F$-{\bf bundle} $\xi\to \mathscr{X}$ over an object $\mathscr{X}$ is a morphism such that 
\begin{enumerate}
\item The fiber admits the structure of $F$.
\item There is an effective epimorphism $p:U\to \mathscr{X}$ such that the pullback of $\xi\to \mathscr{X}$ along $p$ is equivalent to $F\times \mathscr{X}$. 
\end{enumerate}
\end{definition}

The natural home for a global quotient orbifold $\mathscr{X}=M/\!/G$ is in the slice $\infty$-topos $\sh_{\infty}(\cartsp)_{\BB G}$. Indeed, there is a canonical map $M/\!/G\to \ast /\!/G\simeq \BB G$ which projects out $M$, and we can regard this global quotient as an object in the slice. If we take $F$ to be the global quotient given by a $G$-representation $V$, we have a canonical map
$$F:=V/\!/G\to \BB G$$
giving rise to the trivial $V/\!/G$-bundle (or $V$-bundle for short) over the terminal object. Using the general theory of $\infty$-bundles, in particular using the descent axiom for the $\infty$-topos, one finds that a $V$-bundle $\xi\to M/\!/G$ is equivalently the data of a $G$-equivariant vector bundle $\xi^{\prime}\to M$. We can define a $V$-{\bf framing} of $\xi\to M$ in a natural way: namely an equivalence $\phi:\xi\to p_M^*(V/\!/G)$, where $p_M:M/\!/G\to \BB G$ is the canonical map and the superscript denotes the pullback bundle \footnote{Note that the framing depends on the $G$-representation -- there is no canonical choice of $V$.}. We have the following 
\begin{proposition}\label{equitfrm}
A $V$-framing of a bundle $\xi\to M/\!/G$ is equivalently a $G$-equivariant trivialization 
$$\phi:\xi^{\prime}\to V\times M\;,$$
of the associated equivariant bundle $\xi^{\prime}\to M$.
\end{proposition}
\theproof
Consider the diagram in $\sh_{\infty}(\cartsp)$
\(\label{dubpull}
\xymatrix{
V\times M \ar[r]\ar[d] & p_{M}^*(V/\!/G) \ar[r] \ar[d] & V/\!/G\ar[d]
\\
M \ar[r] & M/\!/G\ar[r]^{p_M} & \BB G
}
\)
where $p_M$ projects out $M$ and the left horizontal maps quotient by the $G$-action. The bottom horizontal composite can be written as $M\overset{p_M}{\to}\ast\to \BB G$. Since the fiber of $V/\!/G\to \BB G$ can be canonically identified with $V$, the iterated pullback to $M$ can be identified with $V\times M\to M$. Hence, the outer square in \eqref{dubpull} is cartesian. From the pasting law for pullbacks, the left square is cartesian. By the descent axiom for $\sh_{\infty}(\cartsp)$, it follows also that the top arrow must be the colimiting arrow for a quotient by a $G$-action. Hence $p_M^*(V/\!/G)\simeq (V\times M)/\!/G$. Moreover, this equivalence is canonical in that it is unique, up to a contractible choice. A choice of framing is then a choice of equivalence $\phi:\xi\to (V\times M)/\!/G$. Again, using descent, we identify $\xi\simeq \xi^{\prime}/\!/G$, where $\xi^{\prime}\to M$ is the pullback of $\xi\to M/\!/G$ along $M\to M/\!/G$.  Since the morphism $\phi$ is an arrow in the slice $\sh_{\infty}(\cartsp)_{\BB G}$ by definition, it follows that the data of $\phi$ is equivalent to the data of a $G$-equivariant map $\phi:\xi^{\prime}\to V\times M$. 
\endofproof

Let $\BB{\rm O}(n)$ denote the quotient stack $\ast/\!/{\rm O}(n)$. This stack is the moduli stack of rank $n$-vector bundles with orthogonal structure (see e.g. \cite{SSS}\cite{FSSt}). In particular, for any smooth stack $X$, we have a canonical bijection 
$$
\pi_0\map(X,\BB{\rm O}(n))\cong \pi_0{\rm Vect}(X)\;.
$$
\begin{proposition}
There is a bijective correspondence between $\pi_0\map(M/\!/G,\BB{\rm O}(n))$ and $G$-equivariant vector bundles $V\to M$.
\end{proposition}
\theproof
As discussed above, pullback by maps $f:M/\!/G\to \BB{\rm O}(n)$ gives rise to a bijective correspondence between $\pi_0\map(M/\!/G,\BB{\rm O}(n))$ and isomorphism classes of vector bundles over $M/\!/G$. It remains to show that such isomorphism classes are in canonical bijection with iso classes of equivariant bundles $V\to M$. Let $\xi\to M/\!/G$ be a bundle and let $V\to M$ be the pullback bundle by the quotient $M\to M/\!/G$ and $V$ inherits an action of $G$ by the corresponding action on $M$, giving rise to an equivariant bundle on $M$. this defines a function between iso classes of bundles over $M/\!/G$ and iso classes of equivariant bundles over $M$. To see that this is a bijection, note that by descent, we have an induced bundle equivalence $\xi\to V/\!/G$ over $M/\!/G$. But the induced bundle map $V/\!/G\to M/\!/G$  is just the quotient of the corresponding equivariant bundle $V\to M$. The quotient map by the $G$-action thus induces an inverse function.
\endofproof

The cohomology of global quotient orbifolds $M/\!/G$ is closely related to the equivariant cohomology of the base manifold $M$. Indeed, let
$$\vert \cdot \vert:\sh_{\infty}(\cartsp)\to \mathscr{T}{\rm op}$$
denote the geometric realization functor, sending a smooth stack $X$ to the smooth singular nerve $\vert \mathscr{X} \vert$ (i.e. its geometric realization). Let ${\bf Map}(\mathscr{X},\mathscr{A})$ be the mapping stack between (i.e. the internal hom) smooth stacks $\mathscr{X}$ and $\mathscr{Y}$. In \cite{GH}, it was shown that for an orbifold $\mathscr{X}$ whose isotropy groups are all subgroups of some global group $G$, the geometric realization $\vert {\bf Map}^{f}(\mathscr{X},\mathscr{A})\vert $, with $\mathscr{A}$ an orbifold with isotropy groups in $G$, of the stack of \emph{faithful} (or zero truncated) maps between $\mathscr{X}$ and $\mathscr{A}$ is equivalent to the mapping space between orbi-spaces $\mathscr{P}\mathscr{S}{\rm h}_{\infty}(\mathscr{O}{\rm rb})(y\mathscr{X},y\mathscr{A})$. Here $\mathscr{O}{\rm rb}$ is the global orbit category, with objects the delooping stacks $\BB H$ and morphisms \emph{faithful} maps between them. The functor $y$ is the external Yoneda embedding
$$y:\mathscr{X}\mapsto (\BB H\mapsto \vert {\bf Map}^f(\BB H,\mathscr{X})\vert)\;,$$
where again we take the substack of faithful maps. The faithful mapping stack ${\bf Map}^f(\BB H,\mathscr{X})$ can equivalently be characterized as the \emph{full} mapping stack on the zero truncated morphisms ${\bf Map}_{\BB G}(\mathscr{X},\mathscr{A})$ in the slice topos over the global group $\BB G$ \footnote{For a more detailed discussion on these relationships, we refer the reader to the nlab entry at \newline https://ncatlab.org/nlab/show/orbifold+cohomology.}. From these observations, one sees that the homotopy theory of orbifolds in the slice over $\BB G$ (at least when restricted to the zero truncated objects) is the same as the global equivariant homotopy of fixed point systems prescribed by $y$. If $\mathscr{X}=M/\!/G$ is global quotient orbifold, then a faithful morphism $\BB H\to M/\!/G$ over $\BB G$ simply picks out an $H$-fixed point of $M$. Thus, in this case, we are reduced to the global equivariant homotopy of the $G$-manifold $M$. With these observations, we have the following.

\begin{proposition}\label{orbitequi}
Let $V$ be a $G$ representation and let $S^V=D(V)/S(V)$ be the corresponding representation sphere. Let $H\trianglelefteq G$ be a normal subgroup. We have isomorphisms 
$$\pi_0\vert {\bf Map}_{\BB G}(M/\!/H,S^V/\!/H)\vert \cong \pi^V_{H}(M)\;,$$
where $\pi^V_{H}(M)$ denotes the unstable cohomotopy in ${\rm RO}(G)$-degree \footnote{By $V\in {\rm RO}(G)$ degree, we mean that we are taking equivariant maps to the \emph{representation sphere} $S^V$, with nontrivial action provided by the $H$-representation $V$.} $V$ of the underlying topological space $M$.
\end{proposition}
\theproof
First observe that since $S^V$ is zero truncated, the map $S^V/\!/H\to \BB G$ is zero truncated morphism in the slice. From the main theorem of \cite{GH} and the previous discussion, we have an equivalence
$$\vert {\bf Map}_{\BB G}(M/\!/H,S^V/\!/H)\vert \simeq \vert {\bf Map}^f(M/\!/H,S^V/\!/H)\vert \simeq \mathscr{P}\mathscr{S}{\rm h}_{\infty}(\mathscr{O}{\rm rb})(y(M/\!/H),y(S^V/\!/H))$$
By definition, 
$$y(M/\!/H)\mapsto (\BB K\mapsto \vert {\bf Map}^f(\BB K,M/\!/H)\vert ) \;.$$
By \cite[Lemma 4.4]{GH}, the space $\vert {\bf Map}^f(\BB K,M/\!/H)\vert$ is equivalently the fat realization of the topological groupoid of faithful maps 
$$\left\{\xymatrix{
K\ar@<.05cm>[r] \ar@<-.05cm>[r] & \ast
}
\right\}\to \left\{\xymatrix{
H\times M \ar@<.05cm>[r] \ar@<-.05cm>[r] &M
}
\right\}\;.
$$
This, in turn, is just the degenerate groupoid on the fixed point space $M^K$. The realization is thus equivalent to just $M^K$. If $K\trianglelefteq H$, an element in the mapping space $ \mathscr{P}\mathscr{S}{\rm h}_{\infty}(\mathscr{O}{\rm rb})(y(M/\!/H),y(S^V/\!/H))$ is given by a natural transformation. From the above observation, the corresponding naturality square reduces to 
$$
\xymatrix{
M^K \ar[r]\ar[d] & (S^V)^K\ar[d]
\\
M^L\ar[r] & (S^V)^L
}
$$
where the maps are induced by restriction of fixed points and $K\trianglelefteq L$. From Elemndorf's theorem, we have an equivalence 
$$H\text{-}\mathscr{S}{\rm paces}\overset{\Phi}{\longrightarrow}\mathscr{P}\mathscr{S}{\rm h}_{\infty}(\mathscr{O}{\rm rb}_H)\simeq \mathscr{P}\mathscr{S}{\rm h}_{\infty}(\mathscr{O}{\rm rb})_{\BB H} \;,$$
where the functor $\Phi$ is defined on an $H$-space by $\Phi(M)=(H/K\mapsto M^K)$. By the above observation, we see that the space $ \mathscr{P}\mathscr{S}{\rm h}_{\infty}(\mathscr{O}{\rm rb})(y(M/\!/H),y(S^V/\!/H))$ is equivalent to $\mathscr{P}\mathscr{S}{\rm h}_{\infty}(\mathscr{O}{\rm rb})(u\Phi(M),u\Phi(S^V))$, where $u$ is the fogetful functor $
u:\mathscr{P}\mathscr{S}{\rm h}_{\infty}(\mathscr{O}{\rm rb})_{\BB H}\to \mathscr{P}\mathscr{S}{\rm h}_{\infty}(\mathscr{O}{\rm rb})$. Elmendorf's equivalence then proves the claim.
\endofproof

We introduce the following definition in preparation for our main theorem. 

\begin{definition}
Let $\mathscr{X}=X/\!/H\into M/\!/G=\mathscr{M}$ be a suborbifold and let $V$ be a $G$-representation. A $V$-{\bf framed null-bordism} is an orbifold bordism $\mathscr{W}\into \mathscr{M}\times[0,1]$, equipped with a $V$-framing $\phi$ (of the normal bundle) which restricts to a $V$-framing of its boundary $\mathscr{X}$. In the case $\mathscr{X}=\mathscr{X}\sqcup \mathscr{Y}$, we say that $\mathscr{W}$ is a $V$-framed bordism between $\mathscr{X}$ and $\mathscr{Y}$. We denote the group (under disjoint union) of $V$-framed bordism classes of $k$-dimensional suborbifolds by $\mathscr{F}^k_{V}(\mathscr{M})$ 
\end{definition}

\begin{remark}
Note that in the previous definition, $\mathscr{F}^k_{V}(\mathscr{M})$ is only nonempty for $k={\rm dim}(M)-{\rm dim}(V)$, since a framing is by definition an $H$-equivariant bundle isomorphism $\mathcal{N}\to X\times V$ and in particular $\mathcal{N}_x\cong V$. We will be concerned with the case when $k<{\rm dim}(M)-{\rm dim}(V)$ when we come to our fixed point Pontrjagin-Thom construction and in this case we will weaken the notion of $V$-framing to what we call a partial $V$-framing. Such a framing gives a surjective (but not necessarily injective) bundle map $\mathcal{N}\onto V\times M$
\end{remark}

The previous discussion allows us to study bordisms of orbifolds in terms of equivariant bordisms. In fact, as a consequence of Example \ref{gloqsub}, Propositions \ref{equitorbi} and \ref{equitfrm}, we have the following
\begin{proposition}\label{orbiequi}
Let $H\mathscr{F}^k_{V}(M)$ denote the $H$-equivariant, $V$-framed bordism classes of $k$-dimensional $H$-invariant submanifolds of $M$. We have a bijective correspondence
$$\mathscr{F}^k_{V}(M/\!/G)\cong \bigoplus_{[H]\leq G}H\mathscr{F}^k_{V}(M)\;,$$
where $H$ runs through all conjugacy classes of closed subgroups of $G$.
\end{proposition}
\theproof
For each factor $H\mathscr{F}^k_{V}(M)$, we have a canonical map $H\mathscr{F}^k_{V}(M)\to \mathscr{F}^k_{V}(M)$ which regards a bordism class of $H$-invariant submanifolds as a bordism class of corresponding suborbifolds. This map only depends on the conjugacy class of $H$. Indeed, for $H^{\prime}=gHg^{-1}$ a conjugate, a submanifold $X$ closed under the action of $H^{\prime}$ is ($H^{\prime}$,$H$)-equivariantly diffeomorphic to an $H$-submanifold by the diffeomorphism $g$, i.e. $g(g^{-1}hg)(x)=h(gx)$ for all $x\in X$. Thus the corresponding suborbifolds are isomorphic. Clearly, this map respects the additive structure. This induces a map
$$\bigoplus_{[H]\leq G}H\mathscr{F}^k_{V}(M)\to \mathscr{F}^k_{V}(M/\!/G)\;,$$
which we claim is an isomorphism. But this is clear, since by Proposition \ref{equitorbi}, every suborbifold is a global quotient by an $H$-action for some $H$ and we can take bordisms between these to be $H$-equivariant bordisms. This gives a well defined inverse map 
$$\mathscr{F}^k_{V}(M)\to \bigoplus_{[H]\leq G}H\mathscr{F}^k_{V}(M)\;,$$
which sends each suborbifold to its correspoinding $H$-equivariant bordism class. 
\endofproof

\section{A fixed point Pontrjagin-Thom construction} 

Throughout this section, we fix a compact Lie group $G$ and a compact $G$-manifold $M$. We assume that the quotient Lie groupoid (or action groupoid) $M/\!/G$ is equivalent to an \'etale groupoid. In particular all the isotropy groups $G_x$ are finite and $M/\!/G$ is an orbifold. 

\medskip
We recall the following definition and theorem for $G$-tubular neighborhoods (see e.g. \cite{Ka}).
\begin{definition}
Let $M$ be a $G$-manifold and let $X$ be a $G$-invariant submanifold. A $G$-\emph{invariant tubular neighborhood} of $X$ in $M$ is a pair $(\varphi,\xi)$, where $\xi\to X$ is a $G$-vector bundle over $X$ and $\varphi:\xi\to M$ is a $G$-equivariant embedding onto some open neighborhood $U$ of $X$ in $M$, such that the restriction of $\varphi$ to the zero section of $\xi$ is the inclusion of $X$ into $M$.
\end{definition}

\begin{theorem}[$G$-equivariant tubular neighborhood]
Let $X$ be a closed ($\partial X=\emptyset$) $G$-invariant submanifold of $M$. Then $X$ admits a $G$-invariant tubular neighborhood.
\end{theorem}

Given the previous theorem, we can define the equivariant Pontrjagin-Thom collapse map as follows. Let $H\leq G$ be a subgroup, $V$ be a $G$-representation and let $X\subset M$ be an $V$-framed, $H$-invariant submanifold. Then using the existence of an $H$-equivariant tubular neighborhood $U$ of $X$ in $M$, we can define the collapse map in the usual way, by quotienting out the complement. This gives an $H$-equivariant map 
$${\rm PT}:M\to M/U^c\overset{\varphi}{\simeq} {\rm Th}(\mathcal{N})\overset{\phi}{\simeq} D(V)/S(V)\wedge X_+\overset{{\rm pr}_1}{\to} S^V\;,$$
where $\varphi$ is the map associated with the tubular neighborhood and $\phi$ is induced by the equivariant framing. By the standard differential topology arguments, the homotopy class of this map is independent of a chosen representative for the $G$-bordism class of $X$ and thus gives a well defined map 
\(\label{waserthmabf}
{\rm PT}:\left\{\begin{array}{c}
\text{$H$-Bordism classes of $H$-invariant} 
\\
\text{submanifolds of codimension $d={\rm dim}(V)$}
\\
\text{with $V$-framed normal bundles}
\end{array}\right\}\to \pi^V_H(M)
\)
As observed in \cite{Was}, this map is not an isomorphism in general \footnote{Although this is an iso in the case where $H$ is abelian or finite (see \cite[Theorem 3.11]{Was}).}. In order to justify the statement of our main theorem, we need to better understand why this is the case. 

Let us begin by observing that if $X$ is a $G$-invariant submanifold of $M$, then the slice theorem asserts that for any point $x\in X$, there is a slice $S_x\subset X$ which is invariant under the action of the stabilizer subgroup $G_x$, is transverse to the orbit $Gx$ and $G(S_x)$ is an open neighborhood of the orbit. In essence, the slice theorem tells us how the orbits of $G$ acting on $M$ are stitched together to form a $G$-manifold. The key point is that the Pontrjagin-Thom collapse map knows that these orbits are stitched together smoothly -- the map is defined using only the $G$-manifold $M$ and a choice of $G$-invariant tubular neighborhood. Thus, we cannot expect that if we are given a generic $G$-equivariant map $f:X\to S^V$ that we ought to be able to deform it equivariantly so that $f$ is of the form of ${\rm PT}$. There are simply too few of these types of maps. 

To drive the point home, observe that from the tangent-normal splitting of the embedding $X\into M$, we have $\mathcal{N}\oplus TX\cong TM\vert_{X}$. Since $X$ is $G$-invariant, the existence of a $G$-invariant tubular neighborhood implies that we can choose this decomposition to be compatible with the $G$-action. According to \cite{Was}, given a generic $G$-invariant map $f:M\to S^V$, we can deform $f$ to a map which is in the same homotopy class as ${\rm PT}$, provided that the $G_x$ acts trivially on $T(G/G_x)$ \footnote{This is true if $G$ is finite or abelian.}. The crucial ingredient in proving this claim comes from \cite[Lemma 3.9]{Was}, where we use the fact that in this case the decomposition $T(G/G_x)+S^{\prime}_x=TX_x$ implies that the slice $S^{\prime}_x$ contains the orthogonal complement of the fixed points of $G_x$. This again tells us that $f$ looks similar to a collapse map. Indeed, from the tangent normal decomposition, we have a decomposition $T(G/G_x)+S_x+\mathcal{N}_x=TM_x$, where $S_x$ is the slice in $X$. Though $G$ may act nontrivially on $T(G/G_x)$, $\mathcal{N}_x$ must contain the orthogonal complement of the fixed points, simply because it is orthogonal to $TX_x$ and $X$ is invariant under the action of $G$. Thus, the collapse map will have the property given in \cite{Was}. 
\begin{figure}
\center
\includegraphics[width=.65\textwidth]{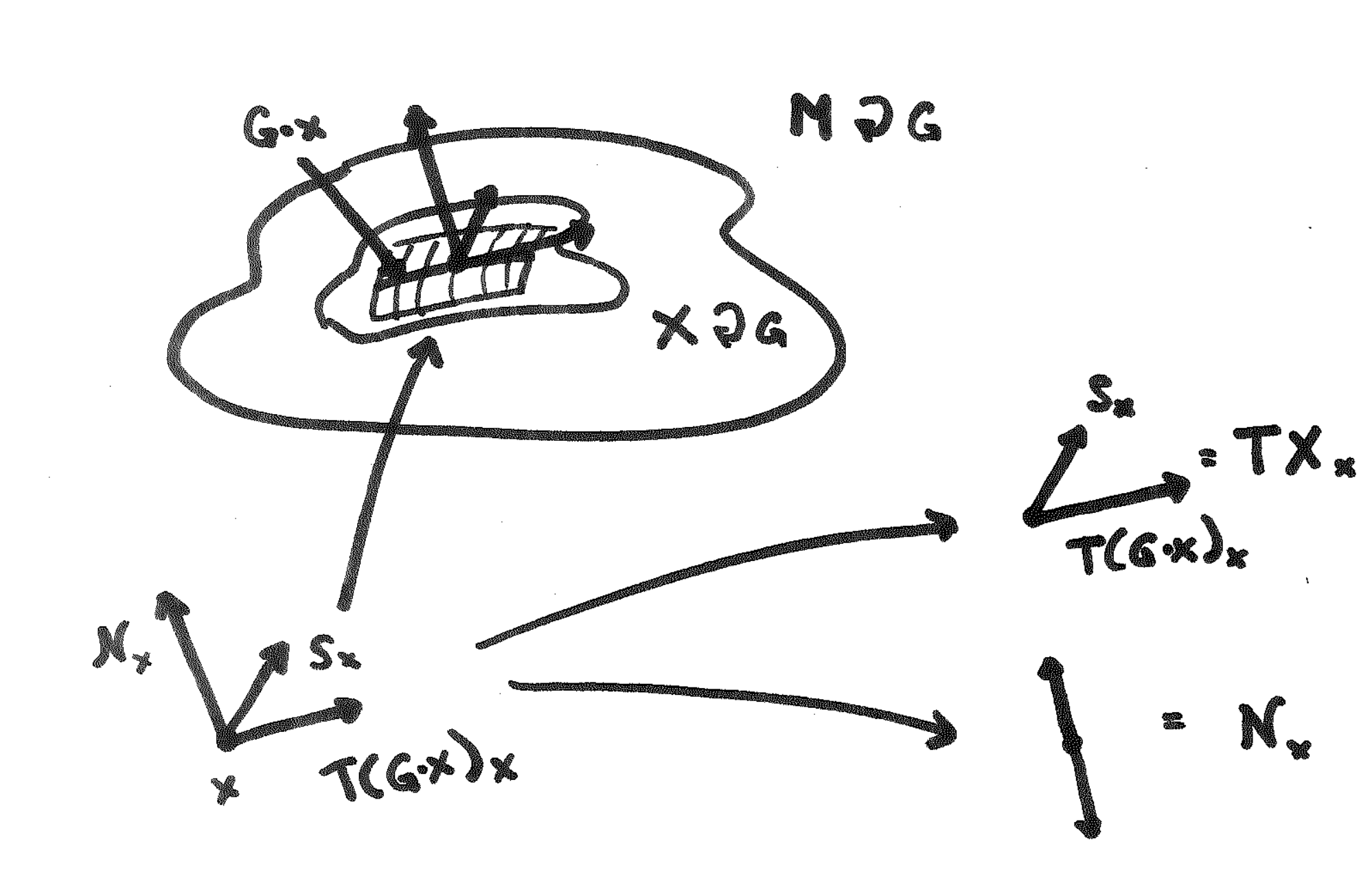}
\caption{Slicing through a $G$-orbit in a $G$-invariant submanifold $X\subset M$}
\end{figure}

The main theorem is proved using the following idea. If we forget that the various orbits of $G$ ought to fit together to form a $G$-invariant manifold (via the slice theorem) then ${\rm PT}$ has a chance to be an isomorphism on the resulting collection of submanifolds. In other words, we focus our attention on each fixed point submanifold. Note however, that if we just consider the fixed point submanifolds and fixed point cobordisms between them, the manifolds will not have the right codimension and there is an obstruction to defining ${\rm PT}$. In order to define the map, we need to modify the usual structure on the normal bundle
\begin{definition}
Suppose $M$ is connected. A {\bf partial $V$-framing} of an equivariant bundle $E\to M$ (of rank possibly different than $V$) is a surjective equivariant bundle map $E\to V$. Equivalently, it is an equivariant decomposition of $E=E_V\oplus E_{T}$ and a choice of $V$-framing of $E_V\to M$. The space of all such lifts is denoted ${\rm Pfr}(E)$ and is identified with the fiber
$$
\xymatrix{
{\rm Pfr}(E)\ar[r]\ar[d] & \ast\ar[d]^-{E}
\\
\map(M/\!/G,\BB{\rm O}(T))\ar[r]^-{i} & \map(M/\!/G,\BB{\rm O}(V\oplus T))
}
$$
where $T$ is an orthogonal complement of $V$ identified as a subspace of $E_x$ and $i:{\rm O}(T)\into {\rm O}(V\oplus T)$ is the canonical map induced by the inclusion as block matrices. For $M=\coprod_
iM_i$ disconnected, a partial $V$-framing of $E\to M$ is a partial $V$-framing of each compenent $E_i\to M_i$
\end{definition}

Note that in the case where $X=X^H$ is an $H$-\emph{fixed point} submanifold of $M$, a partial $V$-framing of an $H$-invariant tubular neighborhood $\mathcal{N}$ of $X$ gives a fiberwise decomposition $\mathcal{N}_x=\mathcal{N}^{\prime}_x+T_x$. In the case where $X$ is a fixed point submanifold of an actual $H$-invariant submanifold $Y$ (of codimension $d={\rm dim}(V))$, the tangent normal decomposition for the embedding $X\into Y$, is compatible with the $G$ action -- as $X$ is fixed. Thus, we have a canonical choice for $T_x$, namely the normal bundle of this embedding. Note that with this choice, $T(G\cdot x)_x\subset T_x$ at $x$, i.e. $\mathcal{N}_T$ contains the distribution defined by the foliation of $X$ by $G$-orbits. One way to see that ${\rm PT}$ cannot be an isomorphism in general is that this distribution may fail to be integrable in general -- i.e. one cannot always find a $V$-framed $G$-invariant submanifold $Y$ (of codimension $d$) containing $X$ as $G$-fixed points. In the present case (i.e. the normal bundle of an embedding) we have the dimension count
$$
{\rm dim}(T) = {\rm dim}(M)-{\rm dim}(X^H)-{\rm dim}(V)\;. \nonumber
$$

\begin{definition}
Let $M$ be an $H$-manifold. A  {\bf partial $V$-framed, fixed point bordism} is a submanifold $W^n\subset M^H$ with boundary the disjoint union of two fixed point submanifolds $M^{n-1},N^{n-1}\subset M^H$, a choice of partial $V$-framing of the normal bundles $\mathcal{N}_{M}$ and $\mathcal{N}_{N}$ (taken in $M$!) and a partial $V$-framing on $\mathcal{N}_{W}$ in $M\times [0,1]$ which restricts to those of $\mathcal{N}_{M}$ and $\mathcal{N}_{N}$. We denote the group of partial $V$-framed, fixed point bordism classes of submanifolds of $M$ by $\mathscr{PF}_{V}(M^H)$. 
\end{definition}

\begin{remark}
\noindent {\bf i}. Given the correspondence between suborbifolds and $H$-invariant submanifolds of $M$, we can define a partial $V$-framing of a suborbifold $\mathscr{X}=X/\!/H$ in a natural way: namely, as a partial $V$-framing of its cover $X$. We denote the group of partial $V$-framed suborbifolds by $\mathscr{PF}_{V}(M/\!/G)$.

\medskip
\noindent {\bf ii}. Note that the subbundle $T\subset \mathcal{N}$ (defined by a partial $V$-framing) and its rank depend on the action of $H$ on $M$. Hence two different subgroups $H,H^{\prime}\leq G$, will correspond to different $T$'s and will have possibly different dimensions (e.g. $H=1$ corresponds to $T=0$). 
\end{remark}

Given a choice of partial $V$-framing, we have a well defined collapse map for \emph{fixed point} submanifolds, given as follows.
\begin{definition}
Let $X^H\into M^H$ be a fixed point submanifold of dimension 
$$d={\rm dim}(M^H)-{\rm dim}(V^H)$$ 
Let $U\overset{\varphi}{\cong} \mathcal{N}\to X^H$ be a $H$-invariant tubular neighborhood of $X^H$ in $M$ and suppose $\mathcal{N}$ admits a partial $V$-framing $\phi$. Then we define ${\rm PT}_{\rm fix}$ by the composite
$${\rm PT}_{\rm fix}:M\to M/U^c\overset{\varphi}{\simeq} {\rm Th}(\mathcal{N})\overset{{\rm Th}(\phi)}{\longrightarrow} D(V)/S(V)\wedge X^H_+\overset{{\rm pr}_1}{\longrightarrow} S^V\;.$$
\end{definition}

Since any homotopy between $H$-fixed point submanifolds can be extended to an $H$-equivariant homotopy (by \cite[Lemma 3.2]{Was}), it follows immediately from the construction that ${\rm PT}_{\rm fix}$ gives a well defined map 
\(\label{fxpbord}
{\rm PT}_{\rm fix}:\mathscr{PF}_{V}(M^H)\to \pi^V_{H}(M)\;,
\)
simply by applying the construction to a fixed point bordism with partial $V$-framing. 

\begin{remark}
There are some subtle dimension counts which could be a source of confusion in the previous definition. For ease of readability, we note that if $X^H\subset (G\acts M)$ denotes the fixed point submanifold with partial $V$-framing $\phi$, then the partial $V$-framing defines a subbundle $T\subset \mathcal{N}_{X^H}$ of the normal bundle $\mathcal{N}_{X^H}$. In this case, we have the dimension count 
\begin{align}
{\rm dim}(X^H) &={\rm dim}(M^H)-{\rm dim}(V^H)\;, \nonumber
\\
{\rm dim}(T) &=({\rm dim}(M)-{\rm dim}(V))-({\rm dim}(M^H)-{\rm dim}(V^H)) \;, \nonumber
\\
{\rm dim}(\mathcal{N}_{X^H}) &= {\rm dim}(T ) +{\rm dim}(V) \;.\nonumber
\end{align}
\end{remark}
Combining this construction with our discussion on orbifolds, we can now prove the main theorem.

\begin{theorem}\label{mainth}
Let $M$ be a smooth, compact $G$-manifold and $H\trianglelefteq G$ a closed normal subgroup. Let $M/\!/G$ denote the resulting global quotient orbifold. We have the following identifications
\begin{enumerate}
\item The fixed point Pontrjagin-Thom construction \eqref{fxpbord} is an isomorphism.
\item We have an isomorphism
$$\pi_0\vert {\bf Map}_{\BB G}(M/\!/H,S^V/\!/H)\vert \cong  \pi_H^V(M)$$
where $\pi^V_H(M)$ is the unstable, $H$-equivariant cohomotopy in ${\rm RO}(H)$-degree $V$.
\item There is an isomorphism 
$$
\bigoplus_{[H]\leq G}\mathscr{PF}_V(M^H)\cong \mathscr{PF}^{\rm fix}_{V}(M/\!/G)\;,$$
where on the right we have partial $V$-framed bordism classes of suborbifolds of the form $M^H\times \BB H$ (i.e. the cover $M^H$ is fixed under the action of $H$).
\end{enumerate}
\end{theorem}
\theproof
Claim 2. is Proposition \ref{orbitequi} and claim 3 follows from Proposition \ref{orbiequi}. Hence, we need only prove the first claim. To do this, we will construct the inverse of ${\rm PT}_{\rm fix}$ \footnote{Note that the the partial $V$-framing is needed to construct ${\rm PT}_{\rm fix}$.}. Fix an $H$-equivariant smooth map $f:M\to S^V$ representing a class in $\pi^V_{H}(M)$. Then restricting $f$ to $M^H$ gives rise to a smooth map $f\vert_{M^H}:M^H\to (S^V)^H$. By the standard transversality argument \cite[Theorem 1.35]{Mil}, $f$ is homotopic to a smooth map $g$ with $0\in (S^V)^H$ as a regular value. By lemma \cite[Lemma 3.3]{Was}, we can extend this homotopy to a $H$-equivariant homotopy $h:f\to g^{\prime}$, such that the restriction of $g^{\prime}$ to $M^H$ agrees with $g$ and has $0$ as a regular value in $(S^V)^{H}$. Let $g^{-1}(0)=X^{H}$ be the corresponding $H$-fixed point submanifold. An equivariant homotopy gives rise to a homotopy on $M^H$ by restriction. Then the standard arguments in differential topology imply that $H$ can be chosen so that $H^{-1}(0)=W\subset M^H$ is a fixed point bordism between the corresponding fixed point submanifolds. Fixing an $H$-invariant tubular neighborhood and restricting $f$ to this neighborhood defines a partial $V$-framing $\phi$. We thus have a well defined homomorphism 
$$\eta:\pi^V_H(M)\to  \mathscr{PF}_{V}(M^H)$$
which sends $\eta:f\mapsto (X_{H},\phi)$. To see that this construction gives rise to a two sided inverse, note that one direction is trivial (i.e. $\eta\circ {\rm PT}_{\rm fix}={\rm id}$). For the other direction, we must show that $\eta$ is injective. Suppose $\eta(f)$ is null-bordant via a bordism with partial $V$-framing. Then we need to show that $f$ is $H$-equivariantly null-homotopic. Application of ${\rm PT}_{\rm fix}$ to this null-bordism produces a null-homotopy for a map which is homotopic to $f$.
\endofproof

In \cite{Was}, the following theorem was proved
\begin{theorem}[Wasserman] 
If $G$ is abelian or finite, then the equivariant Pontrjagin-Thom construction \eqref{waserthmabf} is an isomorphism.
\end{theorem}

As a consistency check we prove the following proposition, which shows that our result is a generalization of Wasserman's. 

\begin{proposition}\label{waschk}
When $H$ is abelian or finite, we have a bijective correspondence
\(\label{fixinvbord}
r:\left\{\begin{array}{c}
\text{Bordism classes of $H$-invariant } 
\\
\text{submanifolds of codim ${\rm dim}(V)$ with }
\\
\text{$V$-framed normal bundles}
\end{array}\right\}\overset{\cong}{\to} 
\left\{\begin{array}{c}
\text{Bordism classes of $H$-fixed point} 
\\
\text{submanifolds of codim ${\rm dim}(V^H)$ in $M^H$ }
\\
\text{ with partial $V$-framed normal bundles}
\end{array}\right\}\;,
\)
where $r$ is the map which restricts from $H$-manifolds to $H$-fixed point submanifolds. The partial $V$-framing is defined by a choice of $H$-equivariant tubular neighborhood of the embedding $X^H\into (H\acts X)$ and a framing of the normal bundle $\mathcal{N}\to X$. 
\end{proposition}
\theproof
First observe that $r$ is well defined since for a bordism $W$ between $H$-manifolds $X$ and $Y$, $W^H$ defines a bordism between $X^H$ and $Y^H$. We need to show that this map is both an epi and monomorphism. Suppose $r(X)=r(Y)$. Then there is a fixed point bordism $Z$ between $X^H$ and $Y^H$ \footnote{We are suppresing the framing information since this just tags along for the ride.}. Applying the fixed point  Pontrjagin-Thom construction gives a corresponding homotopy $H:{\rm PT}_{\rm fix}(X^H)\to {\rm PT}_{\rm fix}(Y^H)$. Since $H$ is abelian or finite, it follows from \cite[Theorem 3.11]{Was} that we can assume $0$ is a regular value of this equivariant homotopy. Then $H^{-1}(0)$ defines an equivariant bordism between $X$ and $Y$. Hence $r$ is injective. To prove surjectivity, observe that by \cite[Theorem 3.11]{Was}, we have a bijective correspondence 
\(\label{wasserbord}
\left\{\begin{array}{c}
\text{Bordism classes of $H$-invariant } 
\\
\text{submanifolds of codim $V$ with }
\\
\text{$V$-framed normal bundles}
\end{array}\right\}\cong \pi^V_{H}(X)
\)
when $H$ is abelian or finite. The inverse map $\eta$ is defined by sending $f:X\to S^V$ to $f^{-1}(0)$ when $f$ is transverse regular to $0\in S^V$. Restricting to fixed points gives a corresponding $H$-fixed point submanifold $f^{-1}(0)^H\subset f^{-1}(0)$. Taking $f={\rm PT}_{\rm fix}$ corresponding to some fixed point submanifold $M^H\subset X^H$ then proves surjectivity.
\endofproof

Relating back to orbifolds, we have the following direct corollary of Proposition \ref{waschk} and Proposition \ref{orbiequi}.
\begin{corollary}
Let $G$ be finite or abelian. Then we have an isomorphism 
$$\mathscr{F}_{V}(G/\!/M)\cong \pi^V_{G}(M)\;,$$
where $\mathscr{F}_{V}(G/\!/M)$ is the group of bordism classes of $V$-framed suborbifolds of the global quotient $M/\!/G$. 
\end{corollary} 

The main theorem was phrased as a correspondence between $H$-fixed bordism classes and the set $\pi^V_H(X)$, with $H$ a subgroup of $G$. Clearly, we could have simply restricted attention to the fixed points of an arbitrary group $G$ and had $H\leq G$ as a special case. However, in the case that $H\trianglelefteq G$ is a normal subgroup, then a fixed point submanifold $X^H$ is $G$-invariant and there are canonical maps
\(\label{indmp}
\mathscr{PF}_{V}(M^H)\to \mathscr{PF}_{V}(M^G)
\)
obtained by taking the $G$-fixed point submanifold of an $H$-fixed point manifold. We also have canonical maps $\pi^V_{G}(M)\to \pi^V_H(M)$ remember only the $H$-equivariance. The two maps go in opposite directions. Note though that from theorem \ref{mainth} we see that the above maps \eqref{indmp} give rise to \emph{induction} maps
$$\pi^V_H(M)\to \pi^V_G(M)$$
on cohomotopy. Asking for a map in the opposite direction as \eqref{indmp} is asking for an \emph{extension} of a $G$-fixed point submanifold to an $H$-fixed point submanifold. In general, we cannot expect such an extension to exist. However, in the case that $G/H$ is finite or abelian, Wasseman's theorem allows us to extend. The fixed point Pontrjagin-Thom construction can be natrually extended to intermediate cases between the two extremes of the full equivariant construction and the fixed point construction. More precisely, for $H\trianglelefteq H^{\prime} \trianglelefteq G$, we can consider the map 
\(\label{intpt}
{\rm PT}^{H}_{H^{\prime}}:M\to M/U^c\overset{\varphi}{\simeq} {\rm Th}(\mathcal{N})\overset{\phi}{\simeq} D(V)/S(V)\wedge X_+\overset{{\rm pr}_1}{\to} S^V\;,
\)
which assigns $H^{\prime}$-equivariant, $H$-fixed point bordisms classes of manifolds, with partial $V$-framing of the normal bundle in $M$.  Let $(H^{\prime}/H)\mathscr{PF}_V(M^H)$ denote the $H^{\prime}$-equivariant bordism classes of $H$-fixed point submanifolds.

\begin{proposition}\label{waserext}
Let $H\trianglelefteq G$ be a closed normal subgroup such that $G/H$ is either finite or abelian. Then there exists extension isomorphisms
$${\rm ext}_{G/H}:\mathscr{PF}_{V}(M^G)\overset{\cong}{\to} (G/H)\mathscr{PF}_{V}(M^H)$$
which extends a $G$-fixed point bordism to a $G$-invariant, $H$-fixed point bordism. Moreover, we have a commutative diagram
\(\label{comdiagwf}
\xymatrix{
\mathscr{PF}_{V}(M^G)\ar[r]^-{{\rm PT}_{\rm fix}}\ar[d]^-{\cong}_{{\rm ext}_{G/H}} & \pi^V_{G}(M)\ar@{=}[d]
\\
(G/H)\mathscr{PF}_{V}(M^H)\ar[r]^-{{\rm PT}^H_{G}} & \pi^V_{G}(M)
}
\)
where the bottom map is defined by the map \eqref{intpt}. 
\end{proposition}
\theproof
Observe that since $G/H$ is abelian, \cite[Theorem 3.11]{Was} implies that we can $G$-homotope a $G$-equivariant map $f:M\to S^V$ to a map which has $0\in S^V$ as a regular value in the $H$-fixed point set $M^H$. This gives a well defined inverse to ${\rm PT}^H_G$ exactly as in the proof of Theorem \ref{mainth}, and ${\rm PT}^H_{G}$ is an isomorphism. It is easy to see from the definition that ${\rm PT}^H_{G}$ agrees with the ${\rm PT}_{\rm fix}$ applied to the corresponding fixed point submanifold $M^G\subset M^H$. Thus the diagram \eqref{comdiagwf} commutes.
\endofproof

Notice that by induction, Proposition \ref{waserext} can be iterated down through a tower of normal subgroups $1\leq H_n\leq H_{n-1}\leq \hdots \leq H_1\leq G$, with abelian (or finite) quotients. In this case, we can construct nested bordism classes of $H_{k-1}/H_k$-invariant submanifolds which extends the bordism class of $G$-fixed point submanifold. That is, we have a commutative diagram 
$$
\xymatrix{
\mathscr{PF}_{V}(M^G)\ar[r]^-{{\rm PT}_{\rm fix}}\ar[d]^-{\cong}_{{\rm ext}_{G/H}} & \pi^V_{G}(M)\ar@{=}[d]
\\
(G/H_1)\mathscr{PF}_{V}(M^{H_1})\ar[r]^-{{\rm PT}^{H_1}_{G}}\ar[d]_{{\rm ext}_{H_1/H_2}} & \pi^V_{G}(M)\ar[d]
\\
(H_1/H_2)\mathscr{PF}_{V}(M^{H_2})\ar[r]^-{{\rm PT}^{H_2}_{H_1}}\ar[d] & \pi^V_{H_1}(M)
\\ 
\vdots\ar[d]_{{\rm ext}_{H_{n-1}/H_n}} &  \vdots \ar[d]
\\
H_{n}\mathscr{F}_{V}(M)\ar[r]^-{{\rm PT}_{H_n}}& \pi^V_{H_n}(M)
}
$$
where the vertical maps on the right are the forgetful maps and at each stage the horizontal maps are isomorphisms. This is reminiscent of one of the fundamental theorems in Galois theory -- namely, solvability of polynomials by radical extensions.

\begin{remark}
 It is interesting to speculate on what the exact analogy here is with Galois theory. There are some hints, for example we can think of fixed point submanifolds as being analogues of field extensions fixed by subgroups of the Galois group.  We can think of submanifolds corresponding to the zero locus of a smooth map $f:M\to S^V$ as the analogue of the zeros of a polynomial. We will not explore here how far this analogy can be stretched.
 \end{remark}

\section{Applications to $M$-theory}

In this section, we work out the example for $G=SU(2)$ acting on the 4-sphere in a way which will be made precise in a moment. The motivation for this example came from the application to $M$-theory \cite{HSS}\cite{FSS1}. We identify the 4-sphere as a representation sphere for $SU(2)$ as follows. Identify $S^{\mathbb{H}}$ as the Thom space of the unit disc in the quaternions $\mathbb{H}$. Via the identification $SU(2)\cong S^3\subset \mathbb{H}$, $SU(2)$ acts by left multiplication on elements of the unit disc and this action stabilizes the boundary $S^3$. Thus $S^{\mathbb{H}}=D(\mathbb{H})/S(\mathbb{H})$ inherits this action. Let $0\in S^{\mathbb{H}}$ denote the image of $0$ under the quotient map $q:D(\mathbb{H})\to D(\mathbb{H})/S(\mathbb{H})$. Then $0$ and its antipodal point are the only fixed point of the $SU(2)$-action -- as $SU(2)$ acts by multiplication on $\mathbb{H}$. For such actions we have the following
\begin{proposition}
Let $H\leq SU(2)$ be a finite subgroup of $SU(2)$ acting on $S^{\mathbb{H}}$ as described above. Then the fixed point Pontrjagin-Thom correspondence gives a bijection
$$
{\rm PT}_{\rm fix}:\left\{\begin{array}{c}\text{Partial ${\mathbb{H}}$-framings of the } 
\\
\text{normal bundle of the full $H$-}
\\
\text{fixed point submanifold $M^H$}
\end{array}\right\}\cong \pi^{\mathbb{H}}_{H}(X)\;.
$$

\end{proposition}
\theproof
Since $H$ only fixes $S^0\in S^\mathbb{H}$, the corresponding fixed point submanifold fills the entirety of the submanifold $M^H$ fixed by $H$ acting on $M$. There is no room for nontrivial bordisms and the statement of theorem \ref{mainth} reduces to the above claim.
\endofproof

It is interesting to see how this statement reflects the physics application. In the above proposition, we should think of the full fixed point submanifold $M^H$ as corresponding to a brane which is fixed by some finite group action $G_{\rm ADE}\subset SU(2)$. For example, in the setting of \cite{HSS}, we consider the MK6-brane living in 11-d space time $M$, identified via the decomposition $M=\RR^{6,1}\times \mathbb{H}$, where we think of the 7-dimensional space $\RR^{6,1}$ (with Lorentzian metric) as the worldvolume of the 6-brane. The submanifold $\RR^{6,1}$ is fixed by various finite group actions $G_{\rm ADE}\subset SU(2)$, acting on $\mathbb{H}$ via left multiplication. In this case, the dimension of the fixed point submanifold happens to have codimension $d={\rm dim}(V)$ in $M$, so the partial $V$-framing just reduces to a $V$-framing of the normal bundle (i.e. an equivariant bundle equivalence $\mathcal{N}\simeq V\times M$). 

\begin{remark}
Note that the spacetime decomposition considered above is non-compact. However, these decompositions are only local, so there is no generality lost in restricting to the compact case.
\end{remark}

\medskip
So far, we have only considered actions coming from the canonical action of $SU(2)$ on $\mathbb{H}$. However, there are other actions which are of physical relevance. Consider the conjugacy action of $SU(2)$ on $\mathbb{H}=\mathbb{R}+\mathbb{H}_{\rm im}$, where $\mathbb{H}_{\rm im}$ denotes the 3-dimensional space of imaginary quaternions. This action is the identity on $\mathbb{R}$ and is the canonical ${\rm SO}(3)$ action on $\mathbb{H}_{\rm im}$ (under the double cover projection ${\rm Spin}(3)\simeq SU(2)\onto {\rm SO}(3)$).  We have a commutative diagram 
$$
\xymatrix{
U(1)\ar@{^{(}->}[r]\ar[d]_-{\cong} & SU(2)\ar[d]
\\
{\rm SO}(2)\ar@{^{(}->}[r] & {\rm SO}(3)
}\;.
$$
A $p$-cyclic subgroup $C_p\subset U(1)$ maps to a cyclic subgroup under the isomorphism on the left. Under the conjugacy action, this cyclic subgroup acts as rotations in a plane in $\mathbb{H}_{\rm im}$. It therefore fixes two directions, the direction perpendicular to the plane and the copy of $\mathbb{R}\into \mathbb{H}$. 

\begin{remark}
Note that in our approach, following \cite{HSS}, we have more than one choice of action. Namely, the action on spacetime and the action on $S^{\mathbb{H}}$, regarded as the target of a field $M\to S^\mathbb{H}$. It is interesting to note that the choice of action on $S^{\mathbb{H}}$ has a major effect on the configuration of the brane in spacetime -- as we indicate in the proceeding discussion.
\end{remark}

From the discussion above, we have the following
\begin{proposition}
The cyclic subgroups $C_p\leq SU(2)$, acting on $\mathbb{H}$ via the conjugacy action give rise to a 2-sphere $S^2\subset S^{\mathbb{H}}$ fixed by the action. In this case, we have a correspondence
$$
{\rm PT}_{\rm fix}:\left\{\begin{array}{c}\text{Bordism classes of $C_p$-fixed point } 
\\
\text{submanifold of codimension 2 in $M^{C_p}$}
\\
\text{with partial ${\mathbb{H}}$-framed normal bundles}
\end{array}\right\}\cong \pi^{\mathbb{H}}_{C_p}(X)
$$
\end{proposition}

The final case we consider is the the dihedral subgroups $D_n\leq SU(2)$. Through the conjugacy action, these groups are generated from cyclic subgroups and one reflection at the plane which these cyclic subgroups act as rotations in. Thus, the orthogonal complement is no longer fixed by the action and the resulting fixed point space is $\mathbb{R}\into \mathbb{H}$. This gives the following. 

\begin{proposition}\label{m5inm6}
The dihedral subgroups $D_n\leq SU(2)$, acting on $\mathbb{H}$ via the conjugacy action give rise to a circle $S^1\subset S^{\mathbb{H}}$ fixed by the action. In this case, we have a correspondence
$$
{\rm PT}_{\rm fix}:\left\{\begin{array}{c}\text{Bordism classes of $D_n$-fixed point } 
\\
\text{submanifold of codimension 1 in $M^{D_n}$}
\\
\text{with partial ${\mathbb{H}}$-framed normal bundles}
\end{array}\right\}\cong \pi^{\mathbb{H}}_{D_n}(X)
$$
\end{proposition}

Returning to the physics applications, this last case is of particular interest. By \cite{AFHS}, the general form of a black M5-brane solution to 11-d supergravity has two limits. In terms of units of the plank length, these are are given as follows:
$$
  \hspace{-.6cm}
  \xymatrix@C=9pt{
    &
    \mbox{
     \footnotesize
         $\left\{\begin{array}{c}\text{
        full}
        \\
        \text{black M5-brane}
        \\
        \text{spacetime}
      \end{array}\right\}$
    }
    \ar[dr]^{ \ \ \  \ell_P  << 1 }
    \ar[dl]_{\ell_P  >> 1\ \ \  }
    \\
    \mathrm{AdS}_7 \times (S^4 /\!/G_{\rm ADE})
    &&
    \mathbb{R}^{5,1} \times \mathrm{C}(S^4 /\!/G_{\rm ADE})
  }
$$
where the left diagonal arrow reflects taking the ``near horizon limit" and the right is the ``far horizon limit". In \cite{AFHS}, it is tacitly assumed that the action on the 4-sphere (in the near horizon limit) is free. However according to \cite[Section 8.3]{MF}, this is not actually the case -- the action of $G_{\rm ADE}$ described there is via the left action of $SU(2)$ acting on $S^{\mathbb{H}}$, which has two fixed points $S^0\subset S^{\mathbb{H}}$. Thus our full fixed point locus in the near horizon limit is $AdS_7\times S^0$, which locally looks like $\RR^{5,1}\times \RR_{>0}\times S^0$. So here the situation is clear -- namely, the action does fix the M-brane, but also fixes a two radial direction spreading out from the removed M5 locus. 

Now in the far horizon limit things become even more pronounced. Here the spacetime geometry changes to that of $\RR^{5,1}\times \mathrm{C}(S^4 /\!/G_{\rm ADE})$, where $\mathrm{C}(S^4 /\!/G_{\rm ADE})$ denotes the metric cone. With the left action through $SU(2)$, we identify $\mathrm{C}(S^4 /\!/G_{\rm ADE})=\RR\oplus \mathbb{H}/\!/G_{\rm ADE}$, with $G_{\rm ADE}$ fixing the subspace $\RR$. Hence, in total, the far horizon spacetime geometry reduces to $\RR^{6,1}\times \HH/\!/G_{\rm ADE}$, which is the MK6-monopole spacetime. 

It is well known (see the literature referenced in Table L of \cite{HSS}) that the M5-brane is a ``domain wall" inside the MK6. This is established mathematically in \cite{HSS}, where it is shown that by intersecting the M5-brane locus with another singularity yields a joint fixed locus which is $\RR^{5,1}$. In terms of the AdS/CFT correspondence, this identifies with the $N=(1,0)$-supersymmetric M5-brane, but not the $N=(2,0)$-supersymmetric M5. 

Going back to Proposition \ref{m5inm6}, we see that there is another way to find the M5-brane inside of the MK6. If, instead of the left action of $SU(2)$ on $S^{\mathbb{H}}$, we take the conjugacy action by a dihedral group, Proposition \ref{m5inm6} says that an $D_n$-equivariant map $M\to S^\mathbb{H}$ will single out  a bordism class of codimension 1-submanifold inside the full fixed point submanifold. According to the above discussion, the \emph{full} fixed point locus can be identified with the MK6. Thus, we find the M5-brane inside the MK6 as this codimension 1 fixed point submanifold. Here though the situation is different, as we did not obtain this 5-dimensional submanifold as the intersection of fixed point loci. Since this is what is ultimately what is responsible for reducing the fermionic dimension by 1/2, we speculate that a supersymmetric counterpart to our main theorem might indeed yield the $N=(2,0)$-supersymmetric M5!

\medskip
It is also interesting to note that very often in these space-time decompositions, there is a radial direction which is singled out -- i.e. a choice of normal direction to the fixed point submanifold. The partial $V$-framing on our fixed point bordisms gives a preferred choice of radial direction, giving also a physical interpretation to this extra structure on the normal bundle. In summary, we have the following correspondences.

\begin{center}
\begin{tabular}{c|c}
{\bf Physics} & {\bf Math}
\\
\hline
brane & fixed point manifold $M^{G_{\rm ADE}}$
\\
\hline
choice of radial direction & choice of distribution $T\subset \mathcal{N}$
\\
\hline
sphere wrapping a brane & unit sphere bundle in $\mathcal{N}=T\oplus \mathcal{N}^{\prime}$
\end{tabular}
\end{center}

\begin{remark}
Note that since the groups here are finite, all the above fixed point submanifolds are contained in an $H$-invariant submanifold of codim 4. Moreover, from proposition \ref{waserext}, the poset of normal subgroups gives rise to a corresponding poset of equivariant bordisms. This poset ought to tells us how various branes are nested (i.e. branes within branes).
\end{remark}

\section*{Acknowledgement}

The author would like to thank John Huerta and Vincent Braunack-Mayer for useful discussions in the initial stages of this project. I am grateful to Hisham Sati and Urs Schreiber for very useful comments on the first draft and for lengthy discussions on the applications to $M$-theory. Finally, I would like to especially thank Urs Schreiber for useful discussions on the relationship between global equivariant cohomology and orbifold cohomology and for providing much of the information on the MK6 and M5 brane in the application section.

\end{document}